\documentclass{article}
\usepackage{latexsym,amsfonts,amsmath,amsthm,makeidx}

\usepackage[left=2.5cm,top=2.5cm,right=2.5cm,bottom=2.5cm]{geometry}

\makeindex
\newtheorem{definition}{Definition}[section]
\newtheorem{theorem}{Theorem}[section]
\newtheorem{lemma}{Lemma}[section]
\newtheorem{corollary}{Corollary}[section]

\newtheorem{remark}{Remark}[section]

\newcommand{\s}{\section}

\newcommand{\R}{\mathbb R}

\newcommand{\bt}{\begin{theorem}}
\newcommand{\et}{\end{theorem}}
\newcommand{\bl}{\begin{lemma}}
\newcommand{\el}{\end{lemma}}
\newcommand{\bd}{\begin{definition}}
\newcommand{\ed}{\end{definition}}
\newcommand{\bc}{\begin{corollary}}
\newcommand{\ec}{\end{corollary}}
\newcommand{\bp}{\begin{proof}}
\newcommand{\ep}{\end{proof}}
\newcommand{\bx}{\begin{example}}
\newcommand{\ex}{\end{example}}
\newcommand{\bi}{\begin{exercise}}
\newcommand{\ei}{\end{exercise}}
\newcommand{\bo}{\begin{prop}}
\newcommand{\eo}{\end{prop}}
\newcommand{\br}{\begin{remark}}
\newcommand{\er}{\end{remark}}
\newcommand{\be}{\begin{equation}}
\newcommand{\ee}{\end{equation}}
\newcommand{\ba}{\begin{align}}
\newcommand{\ea}{\end{align}}
\newcommand{\bn}{\begin{enumerate}}
\newcommand{\en}{\end{enumerate}}
\newcommand{\bg}{\begin{align*}}
\newcommand{\bcs}{\begin{cases}}
\newcommand{\ecs}{\end{cases}}

\newcommand{\NN}{{\mathbb N}}

\newcommand{\bean}{\begin{eqnarray*}}
\newcommand{\eean}{\end{eqnarray*}}

\numberwithin{equation}{section}

\begin{document}

\title{\bf On the equation $p \ \frac{\Gamma(\frac{n}{2}-\frac{s}{p-1})\Gamma(s+\frac{s}{p-1})}{\Gamma(\frac{s}{p-1})\Gamma(\frac{n-2s}{2}-\frac{s}{p-1})}
=\frac{\Gamma(\frac{n+2s}{4})^2}{\Gamma(\frac{n-2s}{4})^2}$}
\date{}
\author{\\{\bf  Senping Luo$^{1}$,\;  Juncheng Wei$^{2}$ \;  and \; Wenming Zou$^{3}$}\\
\footnotesize {\it  $^{1,3}$Department of Mathematical Sciences, Tsinghua University, Beijing 100084, China}\\
\footnotesize {\it  $^{2}$Department of Mathematics, University of British Columbia,}\\
\footnotesize {\it  Vancouver, BC V6T 1Z2, Canada}
}

\maketitle \vskip0.16in
\begin{center}
\begin{minipage}{120mm}
\begin{center}{\bf Abstract}\end{center}
The note is aimed at    giving  a complete characterization of  the following equation:
$$\displaystyle p\frac{\Gamma(\frac{n}{2}-\frac{s}{p-1})\Gamma(s+\frac{s}{p-1})}{\Gamma(\frac{s}{p-1})\Gamma(\frac{n-2s}{2}-\frac{s}{p-1})}
=\frac{\Gamma(\frac{n+2s}{4})^2}{\Gamma(\frac{n-2s}{4})^2}.$$
  The method is based on some key transformation and the properties of the Gamma function. Applications to fractional nonlinear Lane-Emden equations will be given.

\end{minipage}
\end{center}

\vskip0.16in \s{Introduction and main results}
In this note  we consider the following equation on $p$
\be\label{11gamma}
p\frac{\Gamma(\frac{n}{2}-\frac{s}{p-1})\Gamma(s+\frac{s}{p-1})}{\Gamma(\frac{s}{p-1})\Gamma(\frac{n-2s}{2}-\frac{s}{p-1})}
=\frac{\Gamma(\frac{n+2s}{4})^2}{\Gamma(\frac{n-2s}{4})^2},
\ee
where  $p>\frac{n}{n-2s}$ and $ (s, n)$ satisfies
\be
0<s<\frac{n}{2}, \ n\in\NN^+.
\ee

 Equation (\ref{11gamma}) appears frequently    in the study of fractional Lane-Emden equation (see \cite{Wei0=1,Wei1=2}), the fractional Yamabe equation with singularities (see \cite{Del,Gon}) and  also some high-order equations (see \cite{Wei=2,Gazzola2006}, where $s=2$).  For example, consider  the  singular solutions for   the fractional supercritical Lane-Emden equation,
 \begin{equation}
 \label{Lane}
  (-\Delta)^s u =|u|^{p-1}u, \ \ p>\frac{n}{n-2s}.
 \end{equation}
 By Lemma 1.1 of \cite{Wei1=2},   the singular radial solution of (\ref{Lane})  $ u_s$  is given by
  \be
  \label{singular}
  u_s (x)= A |x|^{-\frac{2s}{p-1}}, \ \mbox{where} \ A^{p-1}=\frac{\Gamma(\frac{n}{2}-\frac{s}{p-1})\Gamma(s+\frac{s}{p-1})}{\Gamma(\frac{s}{p-1})\Gamma(\frac{n-2s}{2}-\frac{s}{p-1})}.
  \ee

 By Herbst's generalized Hardy's inequality (\cite{Herb})
 $$ \int_{\R^n} \int_{\R^n} \frac{ (\phi(x)-\phi (y))^2}{|x-y|^{n+2s}} dx dy \geq \frac{\Gamma(\frac{n+2s}{4})^2}{\Gamma(\frac{n-2s}{4})^2} \int_{\R^n} |x|^{-2s} \phi^2 dx,$$
 $u_s$ is stable if and only if the following inequality holds
  \be\label{11gamma1}
p\frac{\Gamma(\frac{n}{2}-\frac{s}{p-1})\Gamma(s+\frac{s}{p-1})}{\Gamma(\frac{s}{p-1})\Gamma(\frac{n-2s}{2}-\frac{s}{p-1})}
 \leq \frac{\Gamma(\frac{n+2s}{4})^2}{\Gamma(\frac{n-2s}{4})^2},
\ee
while it is unstable if
 \be\label{11gamma2}
p\frac{\Gamma(\frac{n}{2}-\frac{s}{p-1})\Gamma(s+\frac{s}{p-1})}{\Gamma(\frac{s}{p-1})\Gamma(\frac{n-2s}{2}-\frac{s}{p-1})}
 > \frac{\Gamma(\frac{n+2s}{4})^2}{\Gamma(\frac{n-2s}{4})^2}.
\ee

In \cite{Wei1=2}, it is proved that  for $ 1<s<2$, if $ p >1, p\not = \frac{n+2s}{n-2s}$ and $p$ satisfies (\ref{11gamma2}) then all stable and finite Morse index solutions to (\ref{Lane}) must be trivial.  An open question is the classification of the range of $p$ for which (\ref{11gamma1}) or (\ref{11gamma2})  holds. In this paper we shall  give an affirmative answer to this question.

Our first result concerns the classification of the roots of (\ref{11gamma}).

\bt\label{11th1} Assume that $n>2s,s>0,p>\frac{n}{n-2s}$.  There exists  $n_0(s) \in \NN^+$ such that
\begin{itemize}
\item [(1)] if  $n\leq n_0(s)$, then \eqref{11gamma} only has one real root  $p$ and  $$\displaystyle p=\frac{n+2s-2+2a_{n,s}\sqrt{n}}{n-2s-2+2a_{n,s}\sqrt{n}},$$
    where $ a_{n,s} $ satisfies  $\frac{1}{\sqrt{n}}<a_{n,s}<\frac{1}{2}\frac{n-2s}{\sqrt{n}}+\frac{1}{\sqrt{n}}$ and is the unique positive root of the function  $f_{n,s} (a)$ defined in (\ref{fnsa}) of Section 2.

\item [(2)] if  $n> n_0(s)$,  then \eqref{11gamma} has  exactly  two real  roots  $p_1$ and $p_2$, where
\be\nonumber
p_1:=\frac{n+2s-2+2a_{n,s}\sqrt{n}}{n-2s-2+2a_{n,s}\sqrt{n}},\;\; \;p_2:=\frac{n+2s-2-2a_{n,s}\sqrt{n}}{n-2s-2-2a_{n,s}\sqrt{n}}.
\ee
Moreover,
 \be\label{roots}
 \frac{n}{n-2s}<p_1<\frac{n+2s}{n-2s}<p_2<+\infty.
 \ee
\end{itemize}
\et

\vskip0.3in

\br The integer $n_0 (s)$ is the largest integer such that
\be
\label{ans}
n-2s -2 - 2 a_{n,s} \sqrt{n} \leq 0.
\ee
Hence when $n>n_0 (s),  n-2s -2 - 2 a_{n,s} \sqrt{n} >0$ and thus $p_2$ is well-defined.
\er

As for the inequalities (\ref{11gamma1})-(\ref{11gamma2}), we have the following sufficient and necessary conditions.

\bt\label{11th2}
Assume that $n>2s,s>0,p>\frac{n}{n-2s}$.  Then there exists $n_0 (s) \in \NN^{+}$ such that
for the inequality \eqref{11gamma2}, we have
\begin{itemize}
\item [(1)]  if  $n\leq n_0(s)$, then the inequality (\ref{11gamma2}) holds  if and only if $$\displaystyle p>p_1:=\frac{n+2s-2+2a_{n,s}\sqrt{n}}{n-2s-2+2a_{n,s}\sqrt{n}};$$

\item [(2)] if  $n> n_0(s)$, then  the inequality (\ref{11gamma2}) holds  if and only if
$$\displaystyle p_1:=\frac{n+2s-2+2a_{n,s}\sqrt{n}}{n-2s-2+2a_{n,s}\sqrt{n}}<p<p_2:=\frac{n+2s-2-2a_{n,s}\sqrt{n}}{n-2s-2-2a_{n,s}\sqrt{n}},$$
where
 \be\nonumber
 \frac{n}{n-2s}<p_1<\frac{n+2s}{n-2s}<\frac{n+2s-4}{n-2s-4}<p_2<+\infty.
 \ee
\end{itemize}
\et

\vskip0.3in

\br\label{June=11} The root $p_1$   appears in the Lane-Emden equation with singularities and also fractional Yamabe equation with singularities, while  the root $p_2$ is  essential   in the study of stability of solutions to the
fractional Lane-Emden equation.  In the literature,  when $s=1$, the root $p_2$ is  usually  called Joseph-Lundgren exponent (See  Joseph and Lundgren \cite{Joseph1972} and  Farina \cite{Farina2007}). When $s=1$, the root $p_1$ plays an important role in constructing singular solutions for Lane-Emden equation with subcritical exponent (Chen-Lin \cite{ChenLin}).

\er

   The following corollary gives a complete classification on the stability of the singular radial solutions to (\ref{Lane})

\bc Assume that $n>2s,s>0,p>\frac{n}{n-2s}$. Let $u_s$ be given by (\ref{singular}). Then there exists  $n_0 (s) \in \NN^{+}$ such that
\begin{itemize}
\item [(1)] if  $n\leq n_0(n,s)$, then  $u_s$ is stable if and only if $ p\geq p_1$;

\item [(2)] if  $n> n_0(n,s)$, then  $u_s$ is stable   if and only if
$$p\geq p_2:=\frac{n+2s-2-2a_{n,s}\sqrt{n}}{n-2s-2-2a_{n,s}\sqrt{n}}.$$

\end{itemize}
\ec

\vskip0.3in

In all the results above, we have the two numbers $n_0(s)$ and $a_{n,s}$ which are to be  implicitly determined.  By (\ref{ans}),  $ n_0(s)$ is the largest integer satisfying the following inequality
\be
\label{n0def}
n \leq  (a_{n,s}+\sqrt{a^2_{n,s}+2s+2})^2
\ee

The range of  $a_{n,s}$ is important in applications. The bound $\frac{1}{\sqrt{n}}<a_{n,s}<\frac{1}{2}\frac{n-2s}{\sqrt{n}}+\frac{1}{\sqrt{n}}$ is too rough. Next, we  give  more
refined and quantitative   estimates on  $a_{n,s}$. These results show that $a_{n,s}$ is very close to the  constant $1$ when $n$ is large.

\bt\label{11th3} Assume that $n>2s,s>0$.
\begin{itemize}
\item [(1)]  For any $\varepsilon_1>0$, there exists $n_1(s,\varepsilon_1)$ such that  $a_{n,s}<1+\varepsilon_1$ whenever $n>\overline n_1(s,\varepsilon_1):=\max\{(1+\varepsilon_1+\sqrt{\max\{(1+\varepsilon_1)^2+2s-2,0\}})^2,n_1(s,\varepsilon_1)\}$,
where $n_1(s,\varepsilon_1)$ is the largest real root of
\be\nonumber\aligned
\Big(
&(-\varepsilon_1^2-2\varepsilon_1)n^4+\big[-27+(18s+48)(1+\varepsilon_1)^2\big]n^3\\
&\;+\big[(-36s^2-96s-144)(1+\varepsilon_1)^2-24s^2-30s+88\big]n^2\\
&\;+\big[(24s^3+192s^2+288s+192)(1+\varepsilon_1)^2+60s^2+64s-144\big]n\\
&\;+48s^4+216s^3+352s^2+288s+192
\Big)=0.
\endaligned
\ee

\item [(2)]  For any $\varepsilon_2>0$, there exists $n_2(s,\varepsilon_2)$ such that $a_{n,s}>1-\varepsilon_2$  whenever  $n>\overline n_2(s,\varepsilon_2):=\max\{(1-\varepsilon_2+\sqrt{\max\{(1-\varepsilon_2)^2+2s-2,0\}})^2,n_2(s,\varepsilon_2)\}$,
where $n_2(s,\varepsilon_2)$ is the square of the largest real root of the following equation $($about variable $t$$)$

\be\nonumber\aligned
&\Big(
(\varepsilon_2^2+2\varepsilon_2)t^6-2(1-\varepsilon_2)^3t^5+\big[18(1-\varepsilon_2)^2-18s-39\big]t^4\\
&+\big[-4(1-\varepsilon_2)^3s-6(1-\varepsilon_2)\big]t^3\\
&+\big[(12s^2+36s)(1-\varepsilon_2)^2+36s^2+144s+158\big]t^2\\
&-12(1-\varepsilon_2)st-24s^3-132s^2-260s-192
\Big)=0.
\endaligned
\ee

\item [(3)]$ \lim_{n\rightarrow+\infty}a_{n,s}=1$\;\;\hbox{for any fixed}\;\;$s>0.$
\end{itemize}
\et

\br  Theorem \ref{11th3} gives   precise  thresholds for $a_{n,s}$.  In fact,  for a fixed range of $s$, say  $s\in (2,3)$, we  have  $0.7<a_{n,s}<1.5$
as long as   $n\geq44$.  Moreover, from the Table 1,  we  have a    quantitative estimate of the constants  $\overline n_1(s,\varepsilon_1)$ and $\overline n_2(s,\varepsilon_2)$ (See Theorem \ref{11th3}).

\begin{table}[!htbp]
\caption{\bf The location  of $a_{n,s}$. For simplicity, we set $a_{n,s}:={\mathcal{A}}.$}
\begin{tabular}{|l|l|l|l|l|l|l|}
\hline
$s \;\&\;  {\mathcal{A}}\;$ & ${\mathcal{A}}<1.5\;$ & ${\mathcal{A}}>0.6\;$ & ${\mathcal{A}}<1.2\;$ & ${\mathcal{A}}>0.8\;$ & ${\mathcal{A}}<1.1\;$ & ${\mathcal{A}}>0.9\;$\\
\hline
$s\in(0,1]$& $n\geq28$ & $n\geq20$     & $n\geq46$ & $n\geq37$ &  $n\geq79$ & $n\geq71$\\
\hline
$s\in(1,2]$ & $n\geq36$ & $n\geq26$    & $n\geq63$ & $n\geq51$ & $n\geq110$ & $n\geq100$\\
\hline
$s\in(2,3]$& $n\geq44$ & $n\geq33$     & $n\geq79$ & $n\geq65$ & $n\geq141$ & $n\geq128$\\
\hline
$s\in(3,4]$ & $n\geq52$ & $n\geq39$    & $n\geq96$ & $n\geq79$ & $n\geq172$ & $n\geq157$\\
\hline
$s\in(4,5]$ & $n\geq59$ & $n\geq46$    & $n\geq112$ & $n\geq93$ & $n\geq204$ & $n\geq186$\\
...&...&...&...&...&...&...
\end{tabular}
\end{table}
\er

\vskip0.3in

\br\label{KKK}  Using the estimates for $a_{n, s}$  we can have some estimates on the critical dimension  $n_0(s)$.  More precisely, for any $\varepsilon_1>0$ we have
\be\label{GGG}
n_0(s)\leq \max\{n_1(s,\varepsilon_1),(1+\varepsilon_1+\sqrt{(1+\varepsilon_1)^2+2s+2})^2\}.
\ee

\vskip0.12in

 If we  select that $\varepsilon_1=1$, we get the Table 2  below.
\er

\begin{table}[!htbp]
\caption{\bf The estimate of $n_0(s)$ for various $s$.}
\begin{tabular}{|l|l|l|l|}
\hline
   $s,n_1(s,\varepsilon_1),n_0(n,s),n_0(s)\quad $ & $n_1(s,\varepsilon_1)<\quad $ & $n_0(n,s)\leq\quad $ &$\quad n_0(s)\leq\quad $ \\
\hline
   $\quad  \quad \quad \quad  s\in(0,1]$& $\quad \quad 22$ & $\quad \quad 24$     & $\quad \quad  24$ \\
\hline
   $\quad \quad \quad \quad  s\in(1,2]$ & $\quad \quad 28$ & $\quad \quad 27$    & $\quad \quad 27$ \\
\hline
   $\quad \quad \quad \quad  s\in(2,3]$& $\quad \quad 33$ & $\quad \quad  30$     & $\quad \quad 33$ \\
\hline
   $\quad \quad \quad \quad  s\in(3,4]$ & $\quad \quad 39$ & $\quad \quad  33$    & $\quad \quad 39$ \\
\hline
  $\quad \quad \quad \quad  s\in(4,5]$ & $\quad \quad 44$ & $\quad \quad 36$    & $\quad \quad 44$ \\
\quad \quad \quad \quad ...&\quad \quad  ...&\quad \quad   ...&\quad \quad ...
\end{tabular}
\end{table}

\newpage

\br  Although the explicit  formula  of $a_{n,s}$ may be very complicated  for general $s>0$,   we have seen  by Theorem   \ref{11th3}-(3)   that
$a_{n,s}$ lies around the constant $1$. Therefore,  the roots $p_1,p_2$ obtained  in Theorem \ref{11th1} and \ref{11th2}
have  the following asymptotic   formulas:
\be\nonumber\aligned
p_1:=\frac{n+2s-2+2a_{n,s}\sqrt{n}}{n-2s-2+2a_{n,s}\sqrt{n}}\approx \frac{n+2s-2+\sqrt{n}}{n-2s-2+\sqrt{n}},\\
p_2:=\frac{n+2s-2-2a_{n,s}\sqrt{n}}{n-2s-2-2a_{n,s}\sqrt{n}}\approx \frac{n+2s-2-\sqrt{n}}{n-2s-2-\sqrt{n}}.
\endaligned\ee
On the other hand, to get more precise estimates on the roots $p_1,p_2$, we just need to  select suitable $\varepsilon_1,\varepsilon_2$ in Theorem \ref{11th3}.

\er

\br  Recall that  when $s=1$, the  Joseph-Lundgren exponent  is given by
the following formula (see Joseph and Lundgren \cite{Joseph1972}, see also Farina \cite{Farina2007} ):
\be\nonumber
p_{JL}(n,s=1):=\begin{cases}
\;\;\;\;\;\;\;\;\;\infty\;\;\;\;\;\;\;\;&\hbox{if}\;  n\leq10,\\
\frac{(n-2)^2-4n+8\sqrt{n-1}}{(n-2)(n-10)}\;\;\;\;\;\;\;\;&\hbox{if}\;  n\geq11.
\end{cases}
\ee
For the bi-harmonic case, i.e.,  $s=2$, Joseph-Lundgren exponent  is given by (see Gazzola and Grunau \cite{Gazzola2006}, see also Davila, Dupaigne, Wang and Wei \cite{Wei=2}):
\be\nonumber
p_{JL}(n,s=2)=\begin{cases}
\;\;\;\;\;\;\;\;\;\;\infty\;\;\;\;\;\;\;\;&\hbox{if}\;\; n\leq12,\\
\frac{n+2-\sqrt{n^2+4-n\sqrt{n^2-8n+32}}}{n-6-\sqrt{n^2+4-n\sqrt{n^2-8n+32}}}\;\;\;\;\;\;\;\;&\hbox{if}\;\; n\geq13.
\end{cases}
\ee
In our setting, we obtain the universal Joseph-Lundgren exponent for any $s>0$, that is,
$$p_{JL}(n,s)=p_2:=\frac{n+2s-2-2a_{n,s}\sqrt{n}}{n-2s-2-2a_{n,s}\sqrt{n}}.$$
In particular,   when $s=1$,
\be
p_{JL}(n,s=1):=\begin{cases}
\;\;\;\;\;\;\;\;\;\infty\;\;\;\;\;\;\;\;&\hbox{if}\;\; n\leq10,\\
\frac{n-2a_{n,1}\sqrt{n}}{n-4-2a_{n,1}\sqrt{n}}\;\;\;\;\;\;\;\;&\hbox{if}\;\; n\geq11,
\end{cases}
\ee
where

\be\nonumber
a_{n,1}=\sqrt{\frac{n-1}{n}}.
\ee
In this case (i.e., $s=1$), the root $p_1$ was obtained in Chen-Lin \cite{ChenLin}  for Lane-Emden equation with subcritical exponent, where
$$p_1=\frac{n+2a_{n,1}\sqrt{n}}{n-4+2a_{n,1}\sqrt{n}}=\frac{n+2\sqrt{n-1}}{n-4+2\sqrt{n-1}}.$$
See Remark \ref{June=11}.

When $s=2$,
\be
p_{JL}(n,s=2)=\begin{cases}
\;\;\;\;\;\;\;\;\;\;\infty\;\;\;\;\;\;\;\;&\hbox{if}\;\;  n\leq12,\\
\frac{n+2-2a_{n,2}\sqrt{n}}{n-6-2a_{n,2}\sqrt{n}}\;\;\;\;\;\;\;\;&\hbox{if}\;\;  n\geq13,
\end{cases}
\ee
where
\be\nonumber
a_{n,2}=\sqrt{\frac{2(n-1)(n^2-2n-2)}{n(n^2+4+n\sqrt{(n-4)^2+4})}}.
\ee
Here we notice that
\be\nonumber
\lim_{n\rightarrow+\infty}a_{n,1}=\lim_{n\rightarrow+\infty}a_{n,2}=1.
\ee

\er

\section{Key transformations and  analysis}
\vskip0.1in

At the first glance,  equation \eqref{11gamma} looks  complicated.  In this section we introduce a key transformation which  puts it in more symmetric form.

First we let
\be\nonumber\frac{2s}{p-1}:=k.\ee
Since $\Gamma(s+1)=s\Gamma(s)$, \eqref{11gamma} becomes
\be\label{11gammak}
\frac{\Gamma(\frac{n-k}{2})\Gamma(s+\frac{k}{2}+1)}{\Gamma(\frac{k}{2}+1)\Gamma(\frac{n-k-2s}{2})}
=\frac{\Gamma(\frac{n+2s}{4})^2}{\Gamma(\frac{n-2s}{4})^2}.
\ee
Here we notice that, the sum of the variables of the Gamma function in both the numerator and  the denominator  on  the left hand side  of the above equation  (\ref{11gammak}) is  equal to$\frac{n+2s}{2}+1$ and  $\frac{n-2s}{2}+1$, respectively.  To make sure that all the variables in the Gamma function in (\ref{11gammak})  have  the term $\frac{1}{4}n+\frac{1}{2}s$ or the term $\frac{1}{4}n-\frac{1}{2}s$, we  introduce a new parameter $a\in \R$ satisfying
\be\label{211}
 k:=\frac{n-(2s+2)}{2}+a\sqrt{n}.\ee
For the reason  that the  term $\sqrt{n}$ appears, see Remark \ref{11rem1}. This is a key point.   Now (\ref{11gammak}) reads as
\be\nonumber\aligned
\frac{\Gamma(\frac{1}{4}n+\frac{1}{2}s+\frac{1}{2}+\frac{1}{2}a\sqrt{n})\Gamma(\frac{1}{4}n+\frac{1}{2}s+\frac{1}{2}-\frac{1}{2}a\sqrt{n})}
{\Gamma(\frac{1}{4}n-\frac{1}{2}s+\frac{1}{2}+\frac{1}{2}a\sqrt{n})\Gamma(\frac{1}{4}n-\frac{1}{2}s+\frac{1}{2}-\frac{1}{2}a\sqrt{n})}
=\frac{\Gamma^2(\frac{1}{4}n+\frac{1}{2}s)}{\Gamma^2(\frac{1}{4}n-\frac{1}{2}s)}.
\endaligned\ee
Now we focus on the new variable $a\in \R$.  Taking the  logarithm on both sides  above we  see that (\ref{11gammak}) becomes
\be\nonumber\aligned
&\underbrace{\ln\Gamma(\frac{1}{4}n+\frac{1}{2}s+\frac{1}{2}+\frac{1}{2}a\sqrt{n})-\ln\Gamma(\frac{1}{4}n+\frac{1}{2}s)}\\
+&\underbrace{\ln\Gamma(\frac{1}{4}n+\frac{1}{2}s+\frac{1}{2}-\frac{1}{2}a\sqrt{n})-\ln\Gamma(\frac{1}{4}n+\frac{1}{2}s))}\\
-&\underbrace{\Big(\ln\Gamma(\frac{1}{4}n-\frac{1}{2}s+\frac{1}{2}+\frac{1}{2}a\sqrt{n})-\ln\Gamma(\frac{1}{4}n-\frac{1}{2}s\Big)}\\
-&\underbrace{\Big(\ln\Gamma(\frac{1}{4}n-\frac{1}{2}s+\frac{1}{2}-\frac{1}{2}a\sqrt{n})-\ln\Gamma(\frac{1}{4}n-\frac{1}{2}s)\Big)}=0.\\
\endaligned\ee
Correspondingly, we denote the left hand side (LHS, for short) of the above equation by the following
\be\nonumber\aligned
LHS:=&g_1(n,s,a)+g_2(n,s,a)-g_3(n,s,a)-g_4(n,s,a)=0;\\
LHS:=&\underbrace{g_1(n,s,a)+g_2(n,s,a)}-\underbrace{(g_3(n,s,a)+g_4(n,s,a))}=0;\\
LHS:=&f_1(n,s,a)-f_2(n,s,a)=0.\\
\endaligned\ee
which can be written as
\be
\label{fnsa}
f_{n,s} (a):= f_1 (n,s,a)- f_2 (n,s, a)=0
\ee

We note  that,  to make sure that all the expressions in the Gamma function above are meaningful, we  need that
$-2<k<n-2s$, equivalently,
\be\label{zou=100}-\frac{1}{2}\frac{n-2s}{\sqrt{n}}-\frac{1}{\sqrt{n}}<a<\frac{1}{2}\frac{n-2s}{\sqrt{n}}+\frac{1}{\sqrt{n}}.\ee
By these notations above, we first observe that $f_{n,s}$ is an even function
\bl\label{11fsym}
\be\nonumber
f_{n,s} (-a)=f_{n, s} (a).
\ee
\el
\bp It can be checked directly. \ep

This nice property allows us to discuss the function $f_{n,s} (a)$ for positive  variable $a\in [0,\frac{1}{2}\frac{n-2s}{\sqrt{n}}+\frac{1}{\sqrt{n}})$.

To obtain further properties of $f_{n, s}(a)$, we introduce the following function
$$ \Psi(x)=\frac{d}{dx}(\ln(\Gamma(x)))=\frac{\Gamma'(x)}{\Gamma(x)}.$$

It is  known that
\be\nonumber
\Psi^{(m)}(x)=(-1)^{m+1}m!\sum_{i=0}^\infty\frac{1}{(x+i)^{m+1}},\;\;m=1,2,...
\ee
For $x>1$, we note that
\be\nonumber\aligned
&\frac{1}{mx^m}=\int_{0}^{+\infty}\frac{1}{(x+y)^{m+1}}dy\leq\sum_{i=0}^\infty\frac{1}{(x+i)^{m+1}}\\
&\leq\int_0^{+\infty}\frac{1}{(x+y-1)^{m+1}}dy=\frac{1}{m(x-1)^m}.
\endaligned\ee
Therefore, by letting $m=2k$ and $m=2k+1$ respectively,  we have the following estimates on the derivatives of $\Psi(x)$:
\be\label{11dpsi}\aligned
\begin{cases}
-\frac{(2k-1)!}{(x-1)^{2k}}\leq \Psi^{(2k)}(x)\leq-\frac{(2k-1)!}{x^{2k}},\;\;m=2k,\\
\frac{(2k)!}{x^{2k+1}}\leq \Psi^{(2k+1)}(x)\leq\frac{(2k)!}{(x-1)^{2k+1}},\;\;m=2k+1.
\end{cases}
\endaligned\ee

\bl\label{11f0} If $n>2s$, $s>0$, then $f_{n,s} (0)>0$.
\el
\bp
Consider the function $\ln\Gamma(\frac{1}{2}n+\frac{1}{4}s+x)-\ln\Gamma(\frac{1}{2}n-\frac{1}{4}s+x)$ for $x\geq0$. We have
\be\nonumber\aligned
&\frac{d}{dx}\Big(\ln\Gamma(\frac{1}{2}n+\frac{1}{4}s+x)-\ln\Gamma(\frac{1}{2}n-\frac{1}{4}s+x)\Big)\\
&=\Psi(\frac{1}{2}n+\frac{1}{4}s+x)-\Psi(\frac{1}{2}n-\frac{1}{4}s+x)\\
&>0
\endaligned\ee
since $s>0$.  It follows that
\be\nonumber\aligned
f_{n, s} (0)=&2\Big(\ln\Gamma(\frac{1}{2}n+\frac{1}{4}s+\frac{1}{2})-\ln\Gamma(\frac{1}{2}n+\frac{1}{4}s)\Big)-\\
&2\Big(\ln\Gamma(\frac{1}{2}n-\frac{1}{4}s+\frac{1}{2})-\ln\Gamma(\frac{1}{2}n-\frac{1}{4}s)\Big)\\
=&2\Big(\ln\Gamma(\frac{1}{2}n+\frac{1}{4}s+x)-\ln\Gamma(\frac{1}{2}n-\frac{1}{4}s+x)\Big)\mid_{x=\frac{1}{2}}\\
-&\;2\Big(\ln\Gamma(\frac{1}{2}n+\frac{1}{4}s+x)-\ln\Gamma(\frac{1}{2}n-\frac{1}{4}s+x)\Big)\mid_{x=0}\\
> &0.
\endaligned\ee
\ep

%=&2\Big(\Psi(\frac{1}{2}n+\frac{1}{4}s+x)-\Psi(\frac{1}{2}n-\frac{1}{4}s+x)\Big)\mid_{x=\frac{1}{2}}\\
%-&\;2\Big(\Psi(\frac{1}{2}n+\frac{1}{4}s+x)-\Psi(\frac{1}{2}n-\frac{1}{4}s+x)\Big)\mid_{x=0}\\

\bl\label{11fdecrease}  Let $a\geq0$. Then $f_{n, s}^{'} (0)=0$ and  $f_{n, s}^{'} (a)<0$ if $a>0$.
\el

\bp

Note that
\be\nonumber\aligned
f_{n,s}^{'} (a)=&\frac{1}{2}\sqrt{n}\Big(\big(\Psi(\frac{1}{4}n+\frac{1}{2}s+\frac{1}{2}a\sqrt{n})-\Psi(\frac{1}{4}n-\frac{1}{2}s+\frac{1}{2}a\sqrt{n})\big)\\
&-\big(\Psi(\frac{1}{4}n+\frac{1}{2}s-\frac{1}{2}a\sqrt{n})-\Psi(\frac{1}{4}n-\frac{1}{2}s-\frac{1}{2}a\sqrt{n})\big)\Big)\\
\endaligned\ee

Since $ f_{n,s} (a)$ is an even function, it follows that $ f_{n, s}^{'} (0)=0$.

For $a>0$, let us consider the function $\Psi(\frac{1}{4}n+x)+\Psi(\frac{1}{4}n-x)$ for $\frac{1}{4}n>x>0$. By (\ref{11dpsi}) we have
\be\nonumber\aligned
\frac{d}{dx}\Big(\Psi(\frac{1}{4}n+x)+\Psi(\frac{1}{4}n-x)\Big)=\Psi'(\frac{1}{4}n+x)-\Psi'(\frac{1}{4}n-x)<0,\;x>0,
\endaligned\ee
then
\be\nonumber\aligned
\frac{d}{da}f_{n, s} (a)=\frac{1}{2}\sqrt{n}\Big(&\big(\Psi(\frac{1}{4}n+x)+\Psi(\frac{1}{4}n-x)\big)\mid_{x=\frac{1}{2}s+a\sqrt{n}}\\
-&\big(\Psi(\frac{1}{4}n+x)+\Psi(\frac{1}{4}n-x)\big)\mid_{x=\frac{1}{2}s-a\sqrt{n}}\Big)<0.
\endaligned\ee

\ep

\bl If $n>2s$, then
\be\nonumber f_{n,s} (a)\mid_{a=\frac{1}{2}\frac{n-2s}{\sqrt{n}}+\frac{1}{\sqrt{n}}}=-\infty.\ee
\el
\bp
If $a=\frac{1}{2}\frac{n-2s}{\sqrt{n}}+\frac{1}{\sqrt{n}}$, by a direct calculation,
 we have that $\frac{1}{4}n-\frac{1}{2}s+\frac{1}{2}-\frac{1}{2}a\sqrt{n}=0$. Thus
 $\ln\Gamma(\frac{1}{4}n-\frac{1}{2}s+\frac{1}{2}-\frac{1}{2}a\sqrt{n})=+\infty$.
Note that
\be\nonumber\aligned
&\ln\Gamma(\frac{1}{4}n-\frac{1}{2}s+\frac{1}{2}+\frac{1}{2}a\sqrt{n})=\ln\Gamma(\frac{1}{2}(n-2s)+1),\\
&\ln\Gamma(\frac{1}{4}n+\frac{1}{2}s+\frac{1}{2}-\frac{1}{2}a\sqrt{n})=\ln\Gamma(s),\\
&\ln\Gamma(\frac{1}{4}n+\frac{1}{2}s+\frac{1}{2}+\frac{1}{2}a\sqrt{n})=\ln\Gamma(\frac{1}{2}n+1).
\endaligned\ee
Therefore,  $f(n,s,a)\mid_{a=\frac{1}{2}\frac{n-2s}{\sqrt{n}}+\frac{1}{\sqrt{n}}}=-\infty$.
\ep
\bc\label{11ccc} If $n>2s$, then
there exists $\varepsilon>0$ small enough such that
\be\nonumber f_{n,s} (a)\mid_{a=\frac{1}{2}\frac{n-2s}{\sqrt{n}}+\frac{1}{\sqrt{n}}-\varepsilon}<0.\ee
\ec
In Lemma \ref{11f0}, we get that $f_{n,s} (0)>0$. Furthermore it holds that
\bl\label{11fsqrtn} If $n>2s$ and $s>0$, then $f_{n,s} (a)\mid_{a=\frac{1}{\sqrt{n}}}>0$.
\el
\bp
\be\nonumber\aligned
&f_{n, s} (a)\mid_{a=\frac{1}{\sqrt{n}}}\\
&=\ln\Gamma(\frac{1}{4}n+\frac{1}{2}s+1)-\ln\Gamma(\frac{1}{4}n+\frac{1}{2}s)\\
&\quad-\Big(\ln\Gamma(\frac{1}{4}n-\frac{1}{2}s+1)-\ln\Gamma(\frac{1}{4}n-\frac{1}{2}s)
\Big).\\
\endaligned\ee
We divide into two different cases.

\medskip

\noindent
Case 1: $s\geq1$. Then we have
\be\nonumber\aligned
&f_{n,s}(\frac{1}{\sqrt{n}})\\
&=\ln\Gamma(\frac{1}{4}n+\frac{1}{2}s+1)-\ln\Gamma(\frac{1}{4}n+\frac{1}{2}s)\\
&\quad-\Big(\ln\Gamma(\frac{1}{4}n-\frac{1}{2}s+1)-\ln\Gamma(\frac{1}{4}n-\frac{1}{2}s)
\Big)\\
&=\Psi(\frac{1}{4}n+\frac{1}{2}s+\theta_1)-\Psi(\frac{1}{4}n-\frac{1}{2}s+\theta_2),
\endaligned\ee
where $\theta_1, \theta_2\in (0,1)$ by the mean value theorem.
Since $(\frac{1}{4}n+\frac{1}{2}s+\theta_1)-(\frac{1}{4}n-\frac{1}{2}s+\theta_2)=s+\theta_1-\theta_2>0$,
and $x>0$, $\Psi'(x)>0$, we obtain the  conclusion.

\medskip

\noindent
Case 2: $0<s<1$. Then we have
\be\nonumber\aligned
&f_{n, s} (\frac{1}{\sqrt{n}})\\
&=\ln\Gamma(\frac{1}{4}n+\frac{1}{2}s+1)-\ln\Gamma(\frac{1}{4}n-\frac{1}{2}s+1)\\
&\quad-\Big(\ln\Gamma(\frac{1}{4}n+\frac{1}{2}s)-\ln\Gamma(\frac{1}{4}n-\frac{1}{2}s)
\Big)\\
&=s\Big(\Psi(\frac{1}{4}n+\frac{\alpha_1}{2}s+1)-\Psi(\frac{1}{4}n+\frac{\alpha_2}{2}s)\Big),
\endaligned\ee
where $\alpha_1,\alpha_2\in (-1,1)$ by mean value theorem.
Since $(\frac{1}{4}n+\frac{\alpha_1}{2}s+1)-(\frac{1}{4}n+\frac{\alpha_2}{2}s)=1+\frac{\alpha_1-\alpha_2}{2}s>1-s>0$
and $x>0$, $\Psi'(x)>0$, we get the  conclusion.
\ep

Recall  that the function $f_{n,s} (a)$ is well-defined if and only if
 \be\label{zou=200}a\in(-\frac{1}{2}\frac{n-2s}{\sqrt{n}}-\frac{1}{\sqrt{n}},\frac{1}{2}\frac{n-2s}{\sqrt{n}}+\frac{1}{\sqrt{n}}).\ee

\bt\label{zou=300} Assume that $n>2s, s>0$ and \eqref{zou=200} holds.
\begin{itemize}
\item [(1)]   The equation $f_{n,s} (a)=0$ of variable $a$ satisfying \eqref{zou=200}  admits precisely  two real roots which are opposite numbers, we denote them as $\pm a_{n,s}$. Moreover,   $\frac{1}{\sqrt{n}}<a_{n,s}<\frac{1}{2}\frac{n-2s}{\sqrt{n}}+\frac{1}{\sqrt{n}}$.

\item [(2)]  The inequality $f_{n, s} (a)>0$ holds if and only if $-a_{n,s}<a<a_{n,s}$.
\end{itemize}

\et
\bp
The proof follows from Corollary \ref{11ccc}, Lemma \ref{11fsqrtn}, Lemma \ref{11fdecrease} and Lemma \ref{11fsym}.
\ep
\vskip0.1in

Now we return to the variable $k$. Recalling that $k=\frac{n-(2s+2)}{2}+a\sqrt{n}$, by Theorem \ref{zou=300} we immediately  have

\bt\label{11thk} Assume that $n>2s, s>0$ and \eqref{zou=200} holds.
\begin{itemize}
\item [(1)]   The equation \eqref{11gammak} of variable $k$ has and only has two real roots, we denote as
$k_1,k_2$. Moreover,    $$k_1:=\frac{n-(2s+2)}{2}-a_{n,s}\sqrt{n}, \quad k_2:=\frac{n-(2s+2)}{2}+a_{n,s}\sqrt{n},$$
where $\frac{1}{\sqrt{n}}<a_{n,s}<\frac{1}{2}\frac{n-2s}{\sqrt{n}}+\frac{1}{\sqrt{n}}$.

\item [(2)]  The inequality  \eqref{11gamma2} holds if and only if $k_1<k<k_2.$
\end{itemize}
\et

\medskip

Now we turn to the original equation \eqref{11gamma} and the corresponding inequality \eqref{11gamma2}.

\vskip0.31in

\noindent{\bf The proofs of Theorems \ref{11th1}-\ref{11th2}.}
Applying Theorem \ref{11thk} above, we get $k_1,k_2$, where $k_1=\frac{n-(2s+2)}{2}-a_{n,s}\sqrt{n}$, $k_2=\frac{n-(2s+2)}{2}+a_{n,s}\sqrt{n}$.
The only difference  between Theorem \ref{11th1}-\ref{11th2} with Theorem \ref{11thk} is that  in Theorem \ref{11th1}-\ref{11th2} $k>0$.
Since $p>1$ in Theorem \ref{11th1}-\ref{11th2}, recalling that $\frac{2s}{p-1}=k$,  we have $k>0$. However in Theorem \ref{11thk}, the region of $k$, that is $-2<k< n-2s$,  is natural from the fact that
 the Gamma function is positive. It  can be easily checked that
$-2<k_1<n-2s$, $\frac{n-2s}{2}<k<n-2s$ since $\frac{1}{\sqrt{n}}<a_{n,s}<\frac{1}{2}\frac{n-2s}{\sqrt{n}}+\frac{1}{\sqrt{n}}$.
Therefore the solution $k_1$ may be non-positive. So we need to divide into several cases, the borderline determined by the following
equation
\be\nonumber
k_1:=\frac{n-(2s+2)}{2}-a_{n,s}\sqrt{n}=0.
\ee
Solving  this, we have either
\be\nonumber
\sqrt{n}=a_{n,s}-\sqrt{a^2_{n,s}+2s+2}  \;\hbox{ or }\;\sqrt{n}=a_{n,s}+\sqrt{a^2_{n,s}+2s+2}.
\ee
Since $a_{n,s}-\sqrt{a^2_{n,s}+2s+2}<0$, we have that $k_1>0$ if and only if $n>(a_{n,s}+\sqrt{a^2_{n,s}+2s+2})^2$.
The rest of the  proofs follow from Theorem \ref{11thk}. \hfill  $\Box$

\section{The location of $a_{n,s}$ and further discussion}

In the section we focus on the constant  $a_{n,s}$, which is crucial in our discussion above, i.e.,  the critical dimension  $n_0(s)$
and the roots of   $p_1$ and $ p_2$ of  \eqref{11gamma}. In the following, we shall give a lower and upper bound of the
function $f_{n,s}(a)$. By these bounds, we can have better estimates for  $a_{n,s}$. A consequence of  the  result  is that
\be\nonumber
\lim_{n\rightarrow+\infty}a_{n,s}=1\;\;\hbox{for any fixed}\;\; s>0.
\ee

\bl\label{11ftaylor2} For $n>2s+4,n>(a+\sqrt{\max\{a^2+2s-2,0\}})^2$, we give an  upper bound of $f_{n,s} (a)$
\be\label{11ftaylor21}\aligned
&f_{n, s} (a)\\
&\leq s\Big(\frac{1}{\frac{1}{4}n-\frac{s}{2}-1}-\frac{\frac{1}{4}+\frac{1}{4}a^2n}{(\frac{1}{4}n+\frac{s}{2})^2}
+\frac{1}{(\frac{1}{4}n-\frac{s}{2}-1)^3}\frac{(\frac{1}{2}a\sqrt{n}+\frac{1}{2})^3}{3}\Big)\\
&=\frac{4s}{3(n-2s-4)^3(n+2s)^2}
\Big\{(-3a^2+3)n^4+2a^3n^{\frac{7}{2}}+\big[-27+(18s+42)a^2\big]n^3\\
&+(8a^3s+6a)n^{\frac{5}{2}}
+\big[(-36s^2-120s-144)a^2-24s^2-30s+86\big]n^2\\
&+(8a^3s^2+24as)n^{\frac{3}{2}}
+\big[(24s^3+168s^2+288s+192)a^2+60s^2+56s-144\big]n\\
&+24as^2n^{\frac{1}{2}}+48s^4+216s^3+344s^2+288s+192
\Big\}
\endaligned\ee
and the lower bound
\be\label{11ftaylor22}\aligned
& f_{n,s}(a)\\
&\geq s
\Big(
\frac{1}{\frac{1}{4}n+\frac{s}{2}}-\frac{\frac{1}{4}+\frac{1}{4}a^2n}{(\frac{1}{4}n-\frac{s}{2}-1)^2}
-\frac{1}{(\frac{1}{4}n-\frac{s}{2}-1)^3}\frac{(\frac{1}{2}a\sqrt{n}-\frac{1}{2})^3}{3}
\Big)\\
&=\frac{4s}{3(n+2s)(n-2s-4)^3}
\Big(
(-3a^2+3)n^3-2a^3n^{\frac{5}{2}}+(18a^2-18s-39)n^2\\
&\quad +(-4a^3s-6a)n^{\frac{3}{2}}+\big[(12s^2+36s)a^2+36s^2+144s+158\big]n\\
&\quad -12asn^{\frac{1}{2}}-24s^3-132s^2-260s-192
\Big).\\
\endaligned\ee
\el
\br\label{11rem1}
 Here we obtain better estimates through  the transform $k=\frac{n-(2s+2)}{2}+a\sqrt{n}$. The term  $\frac{n-(2s+2)}{2}$ seems natural which
guarantees  that all the variables in the Gamma function of the equation \eqref{11gamma} have the part $\frac{1}{4}n+\frac{1}{2}s$ or $\frac{1}{4}n-\frac{1}{2}s$.
\er

\bp If $n>(a+\sqrt{\max\{a^2+2s-2,0\}})^2$ then all the expression in the Gamma function of the function $f_{n, s} (a)$ are positive.

We perform   the Taylor's expansion of the function $g_j(n,s,a)$:

\be\nonumber\aligned
g_1(n,s,a)=&\Psi(\frac{1}{4}n+\frac{1}{2}s)(\frac{1}{2}+\frac{1}{2}a\sqrt{n})+\Psi'(\frac{1}{4}n
+\frac{1}{2}s)\frac{(\frac{1}{2}+\frac{1}{2}a\sqrt{n})^2}{2!}\\
&+\Psi''(\frac{1}{4}n+\frac{\theta_{11}}{2}s)\frac{(\frac{1}{2}+\frac{1}{2}a\sqrt{n})^3}{3!};
\endaligned\ee
\be\nonumber\aligned
g_2(n,s,a)=&\Psi(\frac{1}{4}n+\frac{1}{2}s)(\frac{1}{2}-\frac{1}{2}a\sqrt{n})+\Psi'(\frac{1}{4}n
+\frac{1}{2}s)\frac{(\frac{1}{2}-\frac{1}{2}a\sqrt{n})^2}{2!}\\
&+\Psi''(\frac{1}{4}n+\frac{\theta_{12}}{2}s)\frac{(\frac{1}{2}-\frac{1}{2}a\sqrt{n})^3}{3!};
\endaligned\ee

\be\nonumber\aligned
g_3(n,s,a)=&\Psi(\frac{1}{4}n-\frac{1}{2}s)(\frac{1}{2}+\frac{1}{2}a\sqrt{n})+\Psi'(\frac{1}{4}n
-\frac{1}{2}s)\frac{(\frac{1}{2}+\frac{1}{2}a\sqrt{n})^2}{2!}\\
&+\Psi''(\frac{1}{4}n-\frac{\theta_{21}}{2}s)\frac{(\frac{1}{2}+\frac{1}{2}a\sqrt{n})^3}{3!};
\endaligned\ee
\be\nonumber\aligned
g_4(n,s,a)=&\Psi(\frac{1}{4}n-\frac{1}{2}s)(\frac{1}{2}-\frac{1}{2}a\sqrt{n})+\Psi'(\frac{1}{4}n
-\frac{1}{2}s)\frac{(\frac{1}{2}-\frac{1}{2}a\sqrt{n})^2}{2!}\\
&+\Psi''(\frac{1}{4}n-\frac{\theta_{22}}{2}s)\frac{(\frac{1}{2}-\frac{1}{2}a\sqrt{n})^3}{3!}.
\endaligned\ee
Adding these up and applying the  mean value theorem, we have
\be\nonumber\aligned
f_{n, s} (a)&=g_1(n,s,a)+g_2(n,s,a)-g_3(n,s,a)-g_4(n,s,a)\\
&=\Psi(\frac{1}{4}n+\frac{1}{2}s)-\Psi(\frac{1}{4}n-\frac{1}{2}s)\\
&\quad+\Big(\Psi'(\frac{1}{4}n+\frac{1}{2}s)-\Psi'(\frac{1}{4}n-\frac{1}{2}s)\Big)(\frac{1}{4}+\frac{1}{4}a^2n)\\
&\quad+\Big(\Psi''(\frac{1}{4}n+\frac{\theta_{11}}{2}s)-\Psi''(\frac{1}{4}n-\frac{\theta_{21}}{2}s)\Big)\frac{(\frac{1}{2}+\frac{1}{2}a\sqrt{n})^3}{3!}\\
&\quad+\Big(\Psi''(\frac{1}{4}n+\frac{\theta_{12}}{2}s)-\Psi''(\frac{1}{4}n-\frac{\theta_{22}}{2}s)\Big)\frac{(\frac{1}{2}-\frac{1}{2}a\sqrt{n})^3}{3!}\\
&=s\Big(\Psi'(\frac{1}{4}n+\frac{\alpha_1}{2}s)+\Psi''(\frac{1}{4}n+\frac{\alpha_2}{2}s)(\frac{1}{4}+\frac{1}{4}a^2n)\\
&\quad+\frac{\theta_{11}+\theta_{21}}{2}\Psi'''(\frac{1}{4}n+\frac{\alpha_3}{2}s)\frac{(\frac{1}{2}+\frac{1}{2}a\sqrt{n})^3}{3!}\\
&\quad+\frac{\theta_{21}+\theta_{22}}{2}\Psi'''(\frac{1}{4}n+\frac{\alpha_4}{2}s)\frac{(\frac{1}{2}-\frac{1}{2}a\sqrt{n})^3}{3!}
\Big),
\endaligned\ee
where $\theta_{ij}\in(0,1), \alpha_k\in(-1,1)$.
Now in view of the derivative estimates of $\Psi(x)$ in \eqref{11dpsi}, we get the upper and lower bounds of $f_{n,s} (a)$.
Notice that
\be\nonumber\aligned
&f_{n,s} (a)\leq c_1(n,s,a)\Big((-3a^2+3)n^4+\hbox{lower-order-term}\Big),\\
&f(n,s,a)\geq c_2(n,s,a)\Big((-3a^2+3)n^3+\hbox{lower-order-term}\Big).\\
\endaligned\ee
Let $n=t^2$. Then all the lower order terms are lower order polynomials, that is,
\be\nonumber\aligned
&f_{n,s}(a)\mid_{n=t^2}\leq c_1(n,s,a)\Big((-3a^2+3)t^8+\hbox{lower-order-term}\Big),\\
&f_{n, s} (a)\mid_{n=t^2}\geq c_2(n,s,a)\Big((-3a^2+3)t^6+\hbox{lower-order-term}\Big).\\
\endaligned\ee
Let $a^2>1$, then there exists $t_1=t_1(a)>0$, such that for any $t>t_1$, there holds $(-3a^2+3)t^8+(\hbox{lower-order-term})<0$.
 Hence for any $a>1$, there exists
$t_1=t_1(a)>0$ such that  $f_{n,s} (a)<0$  when $n>t_1^2$. By Lemma \ref{11f0}, we see  $f_{n,s} (0)>0$. Thus,  there is a point $a_0\in(0,a)$ such that $f_{n,s} (a_0)=0$. Further, since $f_{n,s} (a)$ is non-increasing, then $a_0$ is the only real root on interval $[0,+\infty)$.

\vskip0.21in

On the other hand, for any $a>0$ such that $a^2<1$, then there exists $t_2>0$, such that for any $t>t_2$, there holds $(-3a^2+3)n^3+l.o.t>0$.
Therefore, for $\varepsilon_1,\varepsilon_2>0$, there exist $t_1=t_1(\varepsilon_1),t_2=t_1(\varepsilon_2)>0$, such that when $n>\max\{t_1^2,t_2^2\}$, there holds
\be\nonumber\aligned
&f_{n,s} (1+\varepsilon_1)\mid_{n=t^2}\leq c_1(n,s, \varepsilon_1)\Big(-3(\varepsilon_1^2+2\varepsilon_1)t^8+\hbox{lower-order-term}\Big)<0,\\
&f_{n,s}(1-\varepsilon_2)\mid_{n=t^2}\geq c_2(n,s,\varepsilon_2)\Big(3\varepsilon_2(2-\varepsilon_2)t^6+\hbox{lower-order-term}\Big)>0.\\
\endaligned\ee
Thus, there is an $a\in (1-\varepsilon_2,1+\varepsilon_1)$ such that
$f_{n, s} (a)=0$ for $n>\max\{t_1^2,t_2^2\}$. Besides, since $\varepsilon$ is arbitrary small, we get that $a\rightarrow1$ as $n\rightarrow+\infty$.
\ep

By the inverse  transformation of  $k=\frac{2s}{p-1}$ and $k=\frac{n-(2s+2)}{2}+a\sqrt{n}$, a direct consequence of the above lemma is
the following corollary, which complements   Theorem \ref{11th1}.

\bc For any $\varepsilon_1, \varepsilon_2>0$, there exist $a\in(1-\varepsilon_2,1+\varepsilon_1)$ and $n_0=n_0(\varepsilon_1,\varepsilon_2)$, for $n\geq n_0$, the equation \eqref{11gamma}  has and only has two real roots $p_1$ and $p_2$:
\be\nonumber\aligned
p_1:=\frac{n+2s-2+2a\sqrt{n}}{n-2s-2+2a\sqrt{n}},\;\quad p_2:=\frac{n+2s-2-2a\sqrt{n}}{n-2s-2-2a\sqrt{n}}
\endaligned\ee
and $a=a_{n,s}\rightarrow1$ as $n\rightarrow+\infty$.
\ec

\vskip0.12in

\br Generally speaking, when $\varepsilon_1,\varepsilon_2\rightarrow 0+$, $n_0$ will be larger and larger; however,  when $\varepsilon_1,\varepsilon_2$  are far away  from $0$, $n_0$ will be smaller.
That is, to make sure  the existence of such roots, we need to choose the parameters $\varepsilon_1,\varepsilon_2$ suitably  away from $0$.
On the other hand, to get  more accurate estimates on the roots, we need to select $\varepsilon_1,\varepsilon_2\rightarrow 0+$ properly, but that requires   that  $n$ must   be large.
\er

\vskip0.23in
Now we turn to the inequality  \eqref{11gamma2}. By the  transformation above, the inequality \eqref{11gamma2} is equivalent to $f_{n,s}(a)>0$.  Then we have the following
\bc Assume  the inequality \eqref{11gamma2} holds. Then for any $\varepsilon_1, \varepsilon_2>0$, there exists an  $a=a_{n,s}\in(1-\varepsilon_2,1+\varepsilon_1)$ and $n_0=n_0(\varepsilon_1,\varepsilon_2)$, such that for all $n\geq n_0$, we have the following
\be\nonumber\aligned
\frac{n+2s-2+2a\sqrt{n}}{n-2s-2+2a\sqrt{n}}<p<\frac{n+2s-2-2a\sqrt{n}}{n-2s-2-2a\sqrt{n}}
\endaligned\ee
and $a=a_{n,s}\rightarrow1$ as $n\rightarrow+\infty$.
\ec

\br  Generally speaking, when $\varepsilon_1,\varepsilon_2\rightarrow+0$, $n_0$ will be large, while $n_0$ will be small when  $\varepsilon_1,\varepsilon_2$ are  far away from $0$.  Therefore, we must find the balance. In other words,
in order  to get the  existence of such roots, we need to choice the parameter $\varepsilon_1,\varepsilon_2$ suitable away  from $0$;
but to obtain more  accurate estimates of the upper and lower bounds of $a_{n, s}$, we need to select $\varepsilon_1,\varepsilon_2\rightarrow+0$.
\er

\vskip0.2in

To obtain the optimal $n_0(\varepsilon_1,\varepsilon_2)$ and also optimal  upper and lower bound  about $p$ in \eqref{11gamma2},
we  need to generalize Lemma \ref{11ftaylor2}. Again,  we perform the Taylor's expansion of the functions $g_i(n,s,a)$ to $m$ order.

\be\nonumber\aligned
g_1(n,s,a)&=
\sum_{j=0}^m\Psi^{(j)}(\frac{1}{4}n+\frac{1}{2}s)\frac{(\frac{1}{2}+\frac{1}{2}a\sqrt{n})^{j+1}}{{(j+1)}!}\\
&\quad +\Psi^{(m+1)}(\frac{1}{4}n+\frac{\theta_{11}}{2}s)\frac{(\frac{1}{2}+\frac{1}{2}a\sqrt{n})^{m+2}}{(m+2)!};
\endaligned\ee
\be\nonumber\aligned
g_2(n,s,a)&=
\sum_{j=0}^m\Psi^{(j)}(\frac{1}{4}n+\frac{1}{2}s)\frac{(\frac{1}{2}-\frac{1}{2}a\sqrt{n})^{j+1}}{{(j+1)}!}\\
&\quad +\Psi^{(m+1)}(\frac{1}{4}n+\frac{\theta_{12}}{2}s)\frac{(\frac{1}{2}-\frac{1}{2}a\sqrt{n})^{m+2}}{(m+2)!};
\endaligned\ee

\be\nonumber\aligned
g_3(n,s,a)&=
\sum_{j=0}^m\Psi^{(j)}(\frac{1}{4}n-\frac{1}{2}s)\frac{(\frac{1}{2}+\frac{1}{2}a\sqrt{n})^{j+1}}{{(j+1)}!}\\
&\quad+\Psi^{(m+1)}(\frac{1}{4}n-\frac{\theta_{21}}{2}s)\frac{(\frac{1}{2}+\frac{1}{2}a\sqrt{n})^{m+2}}{(m+2)!};
\endaligned\ee
\be\nonumber\aligned
g_4(n,s,a)&=
\sum_{j=0}^m\Psi^{(j)}(\frac{1}{4}n-\frac{1}{2}s)\frac{(\frac{1}{2}-\frac{1}{2}a\sqrt{n})^{j+1}}{{(j+1)}!}\\
&\quad+\Psi^{(m+1)}(\frac{1}{4}n-\frac{\theta_{22}}{2}s)\frac{(\frac{1}{2}-\frac{1}{2}a\sqrt{n})^{m+2}}{(m+2)!}.
\endaligned\ee
Here $\Psi^{(0)}=\Psi$.  Adding these up, then
\be\label{11fnsa}\aligned
&f(n,s,a)\\
&=g_1(n,s,a)+g_2(n,s,a)-g_3(n,s,a)-g_4(n,s,a)\\
&=\sum_{j=0}^m\Big(\Psi^{(j)}(\frac{1}{4}n+\frac{1}{2}s)-\Psi^{(j)}(\frac{1}{4}n-\frac{1}{2}s)\Big)
\frac{(\frac{1}{2}+\frac{1}{2}a\sqrt{n})^{j+1}+(\frac{1}{2}-\frac{1}{2}a\sqrt{n})^{j+1}}{(j+1)!}\\
&\quad+\Big(\Psi^{(m+1)}(\frac{1}{4}n+\frac{\theta_{11}}{2}s)-\Psi^{(m+1)}(\frac{1}{4}n-\frac{\theta_{21}}{2}s)\Big)
\frac{(\frac{1}{2}+\frac{1}{2}a\sqrt{n})^{m+2}}{(m+2)!}\\
&\quad+\Big(\Psi^{(m+1)}(\frac{1}{4}n+\frac{\theta_{12}}{2}s)-\Psi^{(m+1)}(\frac{1}{4}n-\frac{\theta_{22}}{2}s)\Big)
\frac{(\frac{1}{2}-\frac{1}{2}a\sqrt{n})^{m+2}}{(m+2)!}\\
&=s\Big(
\sum_{j=0}^m \Psi^{(j+1)}(\frac{1}{4}n+\frac{\alpha_j}{2}s)
\frac{(\frac{1}{2}+\frac{1}{2}a\sqrt{n})^{j+1}+(\frac{1}{2}-\frac{1}{2}a\sqrt{n})^{j+1}}{(j+1)!}
\Big)\\
&\quad+\frac{\theta_{11}+\theta_{21}}{2}\Psi^{(m+2)}(\frac{1}{4}n+\frac{\alpha_{m+1}}{2}s)\frac{(\frac{1}{2}+\frac{1}{2}a\sqrt{n})^{m+2}}{(m+2)!}\\
&\quad+\frac{\theta_{12}+\theta_{22}}{2}\Psi^{(m+2)}(\frac{1}{4}n+\frac{\alpha_{m+2}}{2}s)\frac{(\frac{1}{2}-\frac{1}{2}a\sqrt{n})^{m+2}}{(m+2)!}
\Big),\\
\endaligned\ee
where $\theta_{ij}\in(0,1),\alpha_k\in (-1,1)$.   When $m=2$, we have the following

\bl\label{11ftaylor3}
Assume  $n>2s+4,n>(a+\sqrt{\max\{a^2+2s-2,0\}})^2$. Then
\be\nonumber\aligned
f_{n,s} (a)&\leq s\Big(
\frac{1}{\frac{1}{4}n-\frac{s}{2}-1}-\frac{\frac{1}{4}+\frac{1}{4}a^2n}{(\frac{1}{4}n+\frac{s}{2})^2}
+\frac{2!}{(\frac{1}{4}n-\frac{s}{2}-1)^3}(\frac{1}{24}+\frac{1}{8}a^2n)
\Big)\\
&=\frac{4s}{3(n-2s-4)^3(n+2s)^2}
\Big(
(-3a^2+3)n^4+\big[-27+(18s+48)a^2\big]n^3\\
&\quad+\big[(-36s^2-96s-144)a^2-24s^2-30s+88\big]n^2\\
&\quad+\big[(24s^3+192s^2+288s+192)a^2+60s^2+64s-144\big] n\\
&\quad+48s^4+216s^3+352s^2+288s+192
\Big).
\endaligned\ee
For the lower bound, we get
\be\nonumber\aligned
f_{n,s} (a)&\geq s\Big(
\frac{1}{\frac{1}{4}n+\frac{s}{2}}-\frac{\frac{1}{4}+\frac{1}{4}a^2n}{(\frac{1}{4}n-\frac{s}{2}-1)^2}
+\frac{2!}{(\frac{1}{4}n+\frac{s}{2})^3}(\frac{1}{24}+\frac{1}{8}a^2n)\\
&\quad-\frac{3!}{(\frac{1}{4}n-\frac{s}{2}-1)^4}(\frac{1}{192}+\frac{1}{192}a^4n^2+\frac{1}{32}a^2n)
\Big)\\
&=\frac{4s}{3(t^2+2s)^3(t^2-2s-4)^4}
\Big(
(-3a^2+3)n^6+\big[-6a^4+(-6s+36)a^2\\
&\quad-12s-51\big] n^5+\big[-36sa^4+(24s^2-276)a^2-12s^2+90s+316\big] n^4\\
&\quad+\big[-72s^2a^4+(48s^3+288s^2+648s+1152)a^2+96s^3+408s^2+\\
&\quad\quad 64s-886\big] n^3\\
&\quad+\big[-48s^3a^4+(-48s^4-768s^3-3312s^2-4608s-3072)a^2-48s^4\\
&\quad-720s^3-2208s^2-1476s+1152\big] n^2
+\big[(-96s^5-192s^4+864s^3\\
&\quad+4608s^2+6144s+3072)a^2-192s^5-816s^4-512s^3\\
&\quad+1656s^2+1536s-1024\big] n
+192s^6+1440s^5+4288s^4\\
&\quad+6224s^3+4608s^2+2048s+1024
\Big)
\endaligned\ee
\el
\br
Combining with the Taylor's expansion of function $f_{n,s}(a)$ in \eqref{11fnsa} and the derivative estimates of $\Psi$ in \eqref{11dpsi}, we can obtain the formula with  higher order expansions. By this way, we can reduce  the bounds $n_1(s,\varepsilon_1)$ and $n_2(s,\varepsilon_2)$.
\er

\noindent{\bf Proof  of Theorems \ref{11th3}.} This is a straightforward consequence of Lemma \ref{11ftaylor2} and Lemma \ref{11ftaylor3}. \hfill$\Box$

\s{\bf Example: an application to $s\in (2,3).$}
\vskip0.1in

Let $p,n,s$ satisfy \eqref{11gamma2}.
From Lemma \ref{11ftaylor2}, we get that
\bl Assume that \eqref{11gamma2} holds.  Then if $n\geq59$, we get that $ 0.7<a<1.5$.  Hence,
\be\label{11LLL}\aligned
\frac{n+2s-2+3\sqrt{n}}{n-2s-2+3\sqrt{n}}<p<\frac{n+2s-2-3\sqrt{n}}{n-2s-2-3\sqrt{n}}
\endaligned\ee
\el

However, if we apply Lemma \ref{11ftaylor3} with higher-order Taylor's expansion, we  may improve the bound $59$ and still get \eqref{11LLL}.
Precisely, we have
\bl Assume that \eqref{11gamma2} holds, then if $n\geq44$, we get that $ 0.7<a<1.5$, hence
\be\nonumber\aligned
\frac{n+2s-2+3\sqrt{n}}{n-2s-2+3\sqrt{n}}<p<\frac{n+2s-2-3\sqrt{n}}{n-2s-2-3\sqrt{n}}.
\endaligned\ee
\el

\end{document}